\documentstyle{amsppt}
\magnification 1200
\UseAMSsymbols
\hsize 6.0 true in
\vsize 8.5 true in
\parskip=\medskipamount
\NoBlackBoxes

\def\mathcal{\Cal}

\def\vp{\varphi}

\def\snint{\raise2pt\hbox{$_{^\not}$}\kern-3.5 pt\int}

\def\snint{\raise2pt\hbox{$_{^\not}$}\kern-3.5 pt\int}

\def\mes{\text {\rm mes\,}}
\TagsOnRight
\NoRunningHeads

\document
\parskip=\medskipamount
\topmatter
\title
On the Lyapounov exponents of Schr\"odinger operators associated with the standard map
\endtitle
\author
J.~Bourgain
\endauthor
\address
\endaddress
\email
\endemail
\abstract
It is shown that Schr\"odinger operators defined from the standard map have positive (mean) Lyapounov exponents for almost all energies
\endabstract
\endtopmatter

\noindent
{\bf 1. Definitions}

(1.1) Let $V$ be a bounded potential on $\Bbb Z_+$.
Define
$$
\align
&\omega(V)=\{W=(W_n)_{n\in\Bbb Z};  \text { there is a sequence } n_j\to\infty \text { such that } d(S^{n_j} V, W)\to 0\}
\\
&\text {where }\\
 &d(V, V')=\sum 2^{-n} |V_n-V_n'|\\
&\text {and }\\
&(S^kV)(n) =V(n+k).
\endalign
$$

\noindent
(1.2) Let $A\subset\Bbb R, \ \mes (A)>0$.
Denote
$$
\Cal R(A)=\{\text{bounded potentials } W=(W_n)_{n\in\Bbb Z} \text { that are reflectionless  on } A\}
$$
and $\Cal R^C(A)=\{W\in\Bbb R(A); \, \sup_n|W_n|\leq C\}$

\noindent
$W\to H=W+\Delta$ and let \ $G(z)$  be the Green's function of $H$. 

Recall that $W$ is reflectionless on $A$ if
$$
 ReG(t)(n)=0 \text {  for a.e. $t\in A$ and every $n\in\Bbb Z$}.
$$
\bigskip

\noindent
{\bf 2.  Remling Theorems} (see \cite {R}).

\noindent
{(2.1)} $\omega(V)\subset \Cal R\big(\sum_{ac}(V)\big)$ \qquad (\cite {R} Theorem 1.4)

\noindent
{(2.2)} $R^C(A)$ is compact \qquad (\cite {R} Proposition 4.1, (d))

\noindent
{(2.3)} The restriction maps $\Cal R(A)\to \Cal R_\pm (A)$ are injective ($\pm$ refers to $\pi_{\Bbb Z_\pm}$).

\noindent
and for any constant $C$, the map
$$
\Cal R^C(A)\to \Cal R^C_+(A) \text { is uniformly continuous}
$$
(\cite {R} Proposition 4.1 (c), (e))

\noindent
(2.4) $W\in \Cal R(A)\Rightarrow A\subset \sum_{ac} (W_\pm)$

\bigskip

\noindent
{\bf 3. Lyapounov Exponents}

Let $T$ be a measure preserving  homeomorphism of a compact metric space $\Omega$ endowed with a probability measure $\mu$ that changes
any non-empty open subset (we do not assume $T$ ergodic). 

Let $\vp\in\Cal C(\Omega)$.

Denote half-line $SO$
$$
H_x=\vp(T^nx)\delta_n+\Delta
$$
and
$$
\align
&M_N(E, x)=\prod^1_N\pmatrix E-\vp (T^nx) &-1\\ 1&0\endpmatrix\\
&L_N(E)=\frac 1N\int\log \Vert M_N(E, x)\Vert\mu (dx)
\endalign
$$
Let $A\subset\Bbb R, |A|>0$ and assume
$$
\underset {N}\to {\underline{\lim}} L_N(E)=0 \text { on } A.\tag 3.1
$$
Then
$$
A\subset\sum\nolimits_{a.c} (V_x) \text { for } \mu \ a.e. x
\tag 3.2
$$

\noindent
{\bf Proof.}

Let
$$
\mu =\int \beta \, \alpha(d\beta)
$$
be the ergodic decomposition of $\mu$.

By Fubini, for $E\in A$
$$
L_N(E)=\int\Big\{\frac 1N\int \log \Vert M_N(E, x)\Vert\beta(dx)\Big\}\alpha (d\beta)
$$
and
$$
\int\Big\{\frac{\lim}N\Big(\frac 1N\int\log \Vert M_N(E, x)\Vert\beta(dx)\Big)\Big\} \alpha (d\beta).
$$
Since $\beta$ is ergodic, it follows that for $\alpha$ - a.e. $\beta$
$$
\frac 1N\log \Vert M_N(E, x)\Vert\to 0 \text { for  $\beta$ - a.e.} \, x\in\Omega
\tag 3.3
$$
Again by Fubini, (3.3) holds for a.e. $E\in A$ and $\beta$ - a.e. \, $x\in\Omega$

By Kotani theory, this implies that
$$
 A\subset\sum\nolimits_{a c} (H_x) \text { for } \beta - a.e. \, x\in\Omega.
\tag 3.4
$$

Since (3.4) is valid for $\alpha$ - a.e. $\beta$, (3.2) follows.
\bigskip

\noindent
{\bf 4. Use of recurrence}

 Let $T$ be as in \S3.

Let $V_x =\big(\vp(T^nx)\big)_{n\in\Bbb Z_+}$

Then
$$
V_x\in\pi_{\Bbb Z_+} \big(\omega(V_x)\big) \text { for  $\mu$  - a.e. } x\in\Omega.\tag 4.1
$$

\noindent
{\bf Proof.}

By Poincar\'e recurrence lemma, for $\mu$ - a.e. $x\in\Omega$, there is a sequence $n_j\to\infty$ such that
$$
T^{n_j}x\to x.\tag 4.2
$$
Hence
$$
\vp(T^{n_j+k}x)\to \vp (T^k x)\equiv W_x(k)
$$
for all $k\in\Bbb Z$ and
$$
d(W_x, V_{T^{n_j}x})\to 0
$$
It follows that
$$
W_x\in\omega(V_x)
$$
and obviously $V_x =\pi_{\Bbb Z_+}(W_x)$.

\bigskip

By (2.1), (4.1) implies
$$
V_x\in \Cal R_+\Big(\sum\nolimits_{a c} (V_x)\Big) \text { for } \mu \ a.e \  x.\tag 4.3
$$
From (4.3), (3.2)
$$
V_x\in \Cal R_+(A) \text { for } \mu \ a.e. \ x\tag 4.4
$$
thus $V_x \in \Cal R^C_+(A)$ for $\mu$  a.e. $x$ and since $\Cal R_+^C(A)$ is compact by (2.2) and $V_x$ depends continuously on $x$, we conclude that
$$
V_x\in\Cal R_+(A) \text { for all } x\in\Omega\tag 4.5
$$
($\Omega$ = closure of any subset of full measure).

Finally, by (2.4)
$$
A\subset \sum\nolimits_{ac}(V_x) \text { for all } x\in\Omega.\tag 4.6
$$

\bigskip

\noindent
{\bf 5. Application to standard map}

Let $T=T_\lambda$ be the standard map on the torus $\Bbb T^2 =\Omega$ with sufficiently large $\lambda$.
Thus 
$$
f_\lambda(x, y) =(-y+2x+\lambda\sin 2\pi x, x).\tag 5.1
$$
Let $\vp\in \Cal C^1(\Bbb T^2)$ be a fixed not-constant function.

If the corresponding SO has vanishing Lyapounov exponents for $E\in A$, (4.6) implies
$$
A\subset\sum\nolimits_{ac} (\{\vp(T^n x)\}_{n\in\Bbb Z_+}\tag 5.2
$$
for all $x\in\Omega$.

By Duarte's work  (theorem A in \cite {D}), there is an invariant hyperbolic set \break 
$\Lambda=\Lambda_\lambda 
\subset\Omega$ such that $T|_\Lambda$ is conjugate to  a Bernoulli shift.
In particular for $x\in \Lambda, (5.2)=\phi$. 
Furthermore, Duarte's result asserts that each point in $\Bbb T^2$ is within a $4\lambda^{-\frac 13}$-neighborhood of $\Lambda$, so that $\vp$ will not be constant on
$\Lambda$ for $\lambda$ large enough.
Hence (5.2) restricted to $\Lambda$ is non-deterministic and Kotani's theorem implies that $A$ is of zero-measure (contradiction).

Hence we proved

\proclaim
{Proposition}
For $\lambda>\lambda_0$, the SO (5.2) associated to the standard map $T_\lambda$ has positive (mean) Lyapounov exponents for a.e. $E\in\Bbb R$.
\endproclaim

\noindent
{\bf Acknowledgement:} The author is grateful to A.~Avila, S.~Sodin and T.~Spencer for some discussions on this topic.

\enddocument